\documentclass[12pt,reqno]{amsart}
\usepackage[notref,notcite]{}
\usepackage{graphicx}
\usepackage{amssymb}
\usepackage{amsmath}

\begin{document}

\def\?{
?\vadjust{\vbox to 0pt{\vss\hbox{\kern\hsize\kern1em\large\bf ?!}}}}

\def\Spin{\mbox{Spin}}
\def\Tr{\mbox{Tr \-}}
\def\inj{\mbox{inj \-}}
\def\R{{\mathbb R}}
\def\C{{\mathbb C}}
\def\H{{\mathbb H}}
\def\Ca{{\mathbb Ca}}
\def\Z{{\mathbb Z}}
\def\N{{\mathbb N}}
\def\Q{{\mathbb Q}}
\def\Ad{\mbox{Ad \-}}
\def\k{{\bf k}}
\def\l{{\bf l}}
\def\sp{\mbox{\bf sp}}

\title[On 4-metrics and evolution]{Evolutions of $S^3$ and $\mathbb{R}P^3$ that describe \protect\\ Eguchi-Hanson metric and metrics of constant curvature}

\author{E.~G.~Malkovich}
\address{Sobolev Institute of Mathematics and Novosibirsk State University,
         Russia}
\email{malkovich@math.nsc.ru}

\thanks{This work was mainly written
while the author was a visiting researcher
at the International Centre for Theoretical Physics (ICTP) in Trieste.
The author is supported by
a~Grant of the Russian Federation for the State Support of Researches
(Contract No.~14.B25.31.0029).}


\begin{abstract}

In this work we illustrate some well-known facts about the evolution of $S^3$ under the
Ricci flow. The Dirac flow we introduce allows us to describe the 4-
dimensional metrics with constant curvature. Another new flow leads to the Eguchi-Hanson metric and can be defined either on metric or on corresponding contact forms.



\end{abstract}

\maketitle
\sloppy

\section{Introduction}

In this article we interpret
some classical 4-dimensional metrics
as deformations of the cone over
$S^3$
(or over
$S^3/\mathbb{Z}_2$)
generated by a certain evolution equations 
of the base.
We define the Dirac flow on a metric on a sphere $S^3$
$$
\frac{\partial}{\partial t} g_{ij}= \sqrt{Ric_{ij}-4Kg_{ij}}, \quad i,j \in \{1,2,3\} \eqno{(DF)}
$$
for $K\in \{-1,~0,~+1\}$ and study it in the simplest situation of conformally round metric. 
If one consider the solution $g=g(t,x)$ of $(DF)$ as a metric on 4-dimensional space then it will be metric of a space with constant curvature.
We also define a~flow on
$\mathbb{R}P^3$, 
which enables us to describe the Eguchi--Hanson metric
$$
\frac{\partial}{\partial t}g_{ij}=\frac{1}{2}{\sqrt{\mathrm{det}(Ric)}}{Ric^{-1}_{ij}}, \quad i,j \in \{ 1,2,3\}.
$$
We show that the evolution of $\mathbb{R}P^3$ described by this flow on metric can also be described by the evolution of contact 1-forms given by the Hitchin flow
$$
(\ast \psi)'=d\psi.
$$
We pose 
the question of constructing new reasonable geometric flows
whose particular solutions would yield
already known metrics with special properties.
Under {\it reasonable} we understand
a~parabolic-type equation on the metric
with a~not-too-complicated dependence of the right-hand side
on the Ricci tensor,
which plays the role of the Laplacian of the metric.
We arrive at this question
pursuing rather naive ideas:
reversing the direction of time,
we can identify the conical singularity
with the singularity developed by the Ricci flow in finite time,
while the existence of asymptotically locally conical metrics
is similar to 
the existence of ancient solutions.

Onda advances similar ideas~\cite{Onda},
asserting in~particular that
one can describe the Taub--NUT metric as
resulting from the action of a~Ricci flow
(or a backward Ricci flow).
He considers the Ricci flow
on various three-dimensional Lie groups
and checks
whether some classical metrics correspond to
solutions to the Ricci flow,
showing no interest in constructing other flows.

\section{The round three-dimensional sphere in
$\mathbb{R}^4$}

To start with,
consider the space
$\mathbb{H}=\mathbb{R}^4$.
The standard flat metric on
$\mathbb{R}^4$
coincides with the cone metric over
$S^3=SU(2)$.
It is known that
the Ricci flow
$\frac{\partial}{\partial t}g(t)=-2 Ric(t)$,
whenever it is defined on the Einstein manifold
$M(t)$,
and on the round sphere in~particular,
changes only the volume of the 
manifold.
Thus,
we can regard
$\mathbb{R}^4$
as the ``configuration space'' of the Ricci flow 
on the round three-dimensional spheres
changing only their radii. Unfortunately the speed of this change is not appropriate, and the metric induced on $\mathbb{R}^4$ will not be flat.
Roughly speaking,
regarding the time coordinate~%
$t$
of some flow
which changes only the radius of the sphere
as the space coordinate~%
$\tau$,
we obtain the space
$(\mathbb{R}_+ \times M(\tau))$
(in this example
$M(t)$
is a sphere of radius~%
$t$)
with strong restrictions on the curvature,
in this case
$R_{ijkl}\equiv0$.

In the case when
$M$
is a~hypersurface
the flow of the mean curvature
$\frac{\partial}{\partial t} g (t)=-2b(t)$
is also defined,
where
$b(t)$
is the second quadratic form depending on the embedding
$r: M \rightarrow H$
of one manifold into the other.
Henceforth we denote the metric and other tensors on the submanifold~%
$M$
of codimension~1
by~%
$\bar{g}$,
$\bar{Ric}$,
and so on,
while tensors on the manifold itself by~%
$g$,
$Ric$,
and so on.
The second quadratic form~%
$b$
is defined only on~%
$M$,
and so we write it without the bar.

Viewing
$S^3$
as the Lie group
$Sp(1)$
of unit quaternions,
we can choose the basis
$(i,~j,~k)$
of three imaginary units
in the tangent space
$T_1 Sp(1)$
at the identity element.
Using the multiplication in 
$Sp(1)$,
extend these 
tangent vectors to three global tangent fields
$(\xi_1,~\xi_2,~\xi_3)$
on the whole sphere.
Associate to them the dual basis in
$T^* S^3=\Lambda^1 (S^3)$
consisting of the 1-forms
$(e^1,~e^2,~e^3)$
with
$e^i(\xi_j)=\delta^i_j$,
usually called the Cartan frame.
Consider on
$\mathbb{R}^4$
the flat conical metric
$$
g=d\tau^2+ \tau^2((e^1)^2+(e^2)^2+(e^3)^2)=d\tau^2+\bar{g}(\tau)
\eqno{(1)}
$$
corresponding to the standard spherical coordinate system,
where~%
$\tau$
plays the role of the radius. 
Firstly,
calculate the second quadratic form of the hypersurface
$S^3\subset \mathbb{R}^4$.

The second quadratic form
$(B)|_{ij}=b_{ij}$
of the embedding
$r: M \rightarrow \mathbb{R}^n$
of a~codimension~1 submanifold
is determined from the equalities
$$
r_{ij}=b_{ij}\mathbf{m}+\Gamma^k_{ij}r_k, \quad i,j,k= 1,\dots n-1
$$
involving the vectors
$r_i=\frac{\partial r}{\partial u^i}$
and
$r_{ij}=\frac{\partial^2 r}{\partial u^i\partial u^j}$
of the first and second partial derivatives of~
$r$,
the unit normal vector~
$\mathbf{m}$,
and the  Christoffel symbols
$\Gamma^k_{ij}$.
Then
$b_{ij}=\langle r_{ij},\mathbf{m}\rangle$.

We can view a~sufficiently small neighborhood of every point of
$r(M)$
as the zero set of some function
$F:\mathbb{R}^n \rightarrow \mathbb{R}$.
Then the gradient of this function
(with respect to the Euclidean metric)
coincides with the normal~%
$\mathbf{m}$
up to a~scalar factor;
therefore,
$\mathbb{R}^n$
includes an~open region diffeomorphic to
$r(M)\times (0,1)$.
Assume that
$\tau=u^0 \in (0,1)$
and consider a~variation
$r_{\tau}$
of the embedding~%
$r$
with 
$\frac{\partial}{\partial \tau}r_{\tau}=r_0=\mathbf{m}$.
Then 
$b_{ij}=\Gamma^0_{ij}$.
Locally these things are rather well-understood;
the main problems arise
when we try to extend the coordinate~%
$\tau$:
the image
$r_{\tau}(M)$
develops degenerations,
which should be controlled.
In our simplest example
$\tau\in (0,\infty)$,
while for
$\tau=0$
the whole sphere collapses to a~point,
the origin.

In the case of an~embedding not into the Euclidean space
$(\mathbb{R}^n,\langle\cdot,\cdot\rangle)$,
but rather into an~arbitrary Riemannian manifold,
we must replace the partial derivatives with covariant derivatives.
Then
$$
\nabla_i\nabla_j r = b_{ij} \mathbf{m} + \Gamma^k_{ij}\nabla_k r.
$$
Observe that
the last formula holds only for the coordinate vector fields
$\nabla_i$
provided that
$[\nabla_i,\nabla_j]=0$,
which fails in our example,
as
$[\xi_i,\xi_{i+1}]=2\xi_{i+2}$,
despite the fact that
the sphere embeds into the Euclidean space~%
$\mathbb{R}^4$.

It is well-known \cite{KoNo}(IV-2) that
the Levi-Civita connection is defined as
$$
\begin{array}{r}
2g(\nabla_X Y,Z)=Xg(Y,Z)+ Yg(X,Z)- Zg(X,Y) +\quad \quad\\ +g([X,Y],Z) +g([Z,X],Y) +g([Z,Y],X).
\end{array}
$$
Putting in our case
$Z=\mathbf{m}=\frac{\partial }{\partial \tau}=\xi_0$,
we obtain
$2B(X,Y)$
in the left-hand side.
Firstly,
calculate
$b_{11}$
for the flat metric~(1).
We have
$$
2b_{11}=2g(\nabla_{\xi_1} \xi_1,\xi_0)=-\xi_0 g(\xi_1,\xi_1)=-2\tau
$$
because only one term in the right-hand side
is not vanishing:
the Lie brackets of fields on
$S^3$
are also fields on
$S^3$
and they remain orthogonal to the radial direction
$\xi_0$,
while the fields
$\xi_i$
are independent of the radial coordinate~%
$\tau$;
thus,
$\xi_i$
commute with
$\xi_0$.
Similarly we calculate two remaining diagonal terms of~
$B$.
The off-diagonal terms vanish identically because
$g(\xi_i,\xi_j)=\delta_{ij}\tau^2$.
It is obvious that
the restriction of the metric~(1) to the three-dimensional sphere of radius~%
$\tau$
satisfies the mean curvature flow:
$$
\Big ( \bar{g}_{ij} \Big )'_{\tau}=\left(
                           \begin{array}{ccc}
                             2\tau & 0 & 0 \\
                             0 & 2\tau & 0 \\
                             0 & 0 & 2\tau \\
                           \end{array}
                         \right)=-2b_{ij}, \quad i,j=1,2,3.
$$

Now 
we would like to check the same equality for the second quadratic form
calculated via the Cartan forms.
Recall that
the connection form of a~Riemannian manifold~%
$M$
is the skew-symmetric matrix
$\omega^i_j$
of 1-forms satisfying
$$
d\varepsilon^i=-\omega^i_j\wedge\varepsilon^j,
$$
where
$\{\varepsilon^1,\ldots,\varepsilon^n\}$
is the orthonormal Cartan (co)frame,
that is,
the basis for
$T^*(M)$
consisting of the 1-forms dual to
the chosen basis of orthonormal tangent vector fields.
In our case
$\{\xi_0,\tau^{-1}\xi_1,\tau^{-1}\xi_2,\tau^{-1}\xi_3\}$
constitute an~orthonormal frame with respect to the metric~%
$g$.
The differentials of the corresponding 1-forms are
$$
\begin{array}{l}
d\varepsilon^0=d(d\tau)=0,\\
d\varepsilon^1=d(\tau e^1)=\frac{1}{\tau}\varepsilon^0\wedge \varepsilon^1+ \frac{2}{\tau}\varepsilon^2\wedge \varepsilon^3,\\
d\varepsilon^2=d(\tau e^2)=\frac{1}{\tau}\varepsilon^0\wedge \varepsilon^2+ \frac{2}{\tau}\varepsilon^3\wedge \varepsilon^1,\\
d\varepsilon^3=d(\tau e^3)=\frac{1}{\tau}\varepsilon^0\wedge \varepsilon^3+ \frac{2}{\tau}\varepsilon^1\wedge \varepsilon^2.\\
\end{array}
$$
Since the matrix
$\omega^i_j$
is skew-symmetric,
we easily calculate that
$$
-(\omega^i_j)|_{i,j=0..3}=\frac{1}{\tau}\left(
                            \begin{array}{cccc}
                              0 & \varepsilon^1 & \varepsilon^2 & \varepsilon^3 \\
                              -\varepsilon^1 & 0 & -\varepsilon^3 & \varepsilon^2 \\
                              -\varepsilon^2 & \varepsilon^3 & 0 & -\varepsilon^1 \\
                              -\varepsilon^3 & -\varepsilon^2 & \varepsilon^1 & 0 \\
                            \end{array}
                          \right).
$$
It is worth noting that
we could also calculate the connection form in the old basis
$\{e^0,e^1,e^2,e^3\}$,
but then we would have to require that
the matrix
$\omega^i_j$,
instead of being skew-symmetric,
that is,
lying in the algebra
$\mathfrak{so}(n)$,
belong to the algebra of matrices
$$
\{ A| AG+GA^\mathrm{T}=0 \},
$$
where~%
$G$
is the matrix of the metric~(1)
not equal to a~scalar matrix.
In other words,
while working with Cartan's structure equations,
we must always use orthonormal frames and coframes.

Since the connection form
is but a~generalization of the Christoffel symbols
$\omega^i_j=\Gamma ^i_{jk}\varepsilon^k$,
everything necessary to calculate the second quadratic form
is already at hand:
$$
b_{jk}=\Gamma^0_{jk}=\omega_j^0(\tau^{-1}\xi_k)=-\frac{1}{\tau}\delta_{jk}, \quad j,k=1,2,3.
$$
Let us verify again,
now in the basis
$\{\varepsilon_1,\varepsilon_2,\varepsilon_3\}$,
that
the standard expanding sphere
satisfies the flow of mean curvature:
$$
 (\bar{g})'_{\tau}=\left ( \varepsilon_1^2+ \varepsilon_2^2+ \varepsilon_3^2 \right )'_{\tau}=\left ( \tau^2 (e_1^2+e_2^2+ e_3^2) \right)'_{\tau}=\frac{2}{\tau}\left ( \varepsilon_1^2+ \varepsilon_2^2+ \varepsilon_3^2 \right)=-2b.
$$
Recall \cite{KoNo}(2-5) that
the curvature form~%
$\Omega$
is the matrix 
of 2-forms
$$
\Omega^i_j =d\omega^i_j +\omega^i_k \wedge \omega^k_j
$$
which generalize the curvature tensor;
namely,
$\Omega ^i_j=\frac{1}{2}R^i_{jkl}\varepsilon^k\wedge \varepsilon^l$.
It is not difficult to verify that~%
$\Omega$
for the metric~(1) vanishes identically;
the space
$\mathbb{R}^4$
under consideration is flat.
Let us calculate the Ricci tensor
$\bar{Ric}$
of~
$\bar{g}$.
By the symmetries of the curvature tensor,
$$
\bar{Ric}_{11}=\bar{R}^1_{111}+\bar{R}^2_{121}+\bar{R}^3_{131}=\bar{R}^1_{212}+ \bar{R}^1_{313}.
$$
In order to calculate the last two terms,
we have to calculate the curvature form restricted to the sphere
$$
\begin{array}{l}
\bar{\Omega}^1_2=d\omega^1_2 +\omega^1_k\wedge \omega^k_2=d(\frac{1}{\tau}\varepsilon^3)+ (-\frac{1}{\tau}\varepsilon^2)\wedge(-\frac{1}{\tau}\varepsilon^1)=
-\frac{1}{\tau^2}\varepsilon^0\wedge \varepsilon^3 \\+ \frac{1}{\tau}(\frac{1}{\tau}\varepsilon^0\wedge\varepsilon^3 +\frac{2}{\tau}\varepsilon^1\wedge\varepsilon^2) +\frac{1}{\tau^2}\varepsilon^2\wedge\varepsilon^1=\frac{1}{\tau^2}\varepsilon^1\wedge\varepsilon^2= \frac{1}{2}\bar{R}^1_{212}\varepsilon^1\wedge\varepsilon^2
\end{array}
$$
with summation over
$k\in\{1,2,3\}$.
Then the components of the Ricci tensor are
$\bar{Ric}_{11}=\frac{4}{\tau^2}=\bar{Ric}_{22}=\bar{Ric}_{33}$.
It may seem here that
we arrive at a~contradiction
since it is known
(see for instance \cite{Be}, 1.159)
that,
when the metric is multiplied by a~constant~%
$\lambda$,
the Riemann (4,0)-tensor does too,
the scalar curvature is multiplied by
$\lambda^{-1}$,
while the Ricci tensor stays unchanged.
All our spheres result from a~sphere of unit radius
by simple homothety,
and
$\lambda\equiv\tau$.
Nevertheless,
we find that
the Ricci tensor depends on~%
$\tau$.
The reason is that
the basis
$\{\varepsilon_1,\varepsilon_2,\varepsilon_3\}$
used to calculate the components
$\bar{Ric}_{ij}$
also depends on~%
$\tau$.
Now we note that in the chosen coordinates
the left-hand side of the Ricci flow equation
$\bar{g}'_{\tau}=-2\bar{Ric}$
is a~positive definite form,
while the right-hand side is negative definite.
Therefore,
as the space coordinate~%
$\tau$
increases,
the metric must ``decrease''
and the corresponding sphere must collapse. 
Consequently,
the standard collapse of the sphere
described by the spherical coordinate system
fails to constitute a~solution to the Ricci flow.

\section[]{
A~metric with one functional parameter}

Consider now the metric 
$$
g=dt^2+f(t)^2\left( \sum_{i=1}^3 (e^i)^2\right )=dt^2 +\bar{g}(t)=\left( \sum_{i=0}^3 (\varepsilon^i)^2\right ). \eqno{(2)}
$$
Here we have $\bar{g}=f(t)^2g_0$, where $f$ is a conformal factor and $g_0$ is a standard metric on $S^3$ that does not depend on $t$ anyhow.
Although now the 1-forms
$\varepsilon^i, ~ i=1,2,3$
depend on a~function~%
$f$,
we keep the previous notation,
being forced to work with orthonormal frames.
Then in the basis
$\{\varepsilon^0,\varepsilon^1,\varepsilon^2,\varepsilon^3\}$
the connection form of the metric~(2) becomes
$$
-\omega=\frac{1}{f}\left(  \begin{array}{cccc}
                              0 & f'\varepsilon^1 & f'\varepsilon^2 & f'\varepsilon^3 \\
                              -f'\varepsilon^1 & 0 & -\varepsilon^3 & \varepsilon^2 \\
                              -f'\varepsilon^2 & \varepsilon^3 & 0 & -\varepsilon^1 \\
                              -f'\varepsilon^3 & -\varepsilon^2 & \varepsilon^1 & 0 \\
                            \end{array}
\right).
$$
An~elementary calculation of the curvature form~%
$\bar{\Omega}$
as above
yields
$\bar{Ric}_{ij}=\frac{4}{f^2}\delta_{ij}$.
To satisfy the Ricci flow equation,
we should require that
$$
\Big ( \bar{g}_{ij} \Big )'_t=\frac{2f'}{f}\delta_{ij}=-\frac{8}{f^2}\delta_{ij}=-2\bar{Ric}_{ij}, \quad \quad i,j=1,2,3;\eqno{(3)}
$$
hence,
$f(t)=\sqrt{8(t_0-t)}$.
This shows that
the radius~%
$f$
of the sphere
must depend on time~%
$t$
as the square root,
and under the action of the Ricci flow
the sphere collapses in finite time,
which is a~well-known fact.
Observe also that
for the chosen basis
the equation of the mean curvature flow
is automatically satisfied;
this is the first equality in~(3).
Therefore,
the mean curvature flow
in fact becomes a~tautology
when we work in an~invariant basis.

Let us calculate the curvature form of the metric~(2)
for the ambient space.
Now
$t$
plays the role of the space coordinate:
$$
\Omega^0_1=-\frac{f''}{f}\varepsilon^0\wedge\varepsilon^1=\frac{1}{2}R^0_{101}\varepsilon^0\wedge\varepsilon^1,
$$
$$
\Omega^1_2=\frac{1-f'^2}{f^2}\varepsilon^1\wedge\varepsilon^2=\frac{1}{2}R^1_{212}\varepsilon^1\wedge\varepsilon^2.
$$
Accordingly,
the Ricci tensor is
$$
Ric_{00}=-6\frac{f''}{f},\quad Ric_{11}=\frac{1}{f^2}(4-4f'^2-2f''f)=Ric_{22}=Ric_{33},
$$
and the scalar curvature is
$$
R=\frac{3}{f^2}(4-4f'^2-4f''f).
$$
For
$f(t)=f_0(t)=t$
we obtain the flat space
$\mathbb{R}^4$,
for
$f(t)=f_{+1}(t)=\sin(t)$
the round sphere
$S^4$
of radius~1,
and for
$f(t)=f_{-1}(t)=\sinh(t)$
the hyperbolic space
$H^4$.

However,
when we consider the Ricci flow with respect to time~%
$\tau$,
for an~appropriately chosen dependence between 
$t$~%
and~
$\tau$
the sphere
$S^3$
under the action of this flow
sweeps a~space of constant curvature.
For the flat space
$\mathbb{R}^4$
we must put
$\frac{dt}{d\tau}=-\frac{4}{t}$
or
$\tau=h_0=const-\frac{t^2}{8}$.
For the sphere
$S^4$
we must put
$\tau=h_{+1}=const+\frac{1}{8}\sin^2{t}$,
and for the hyperbolic space
$\tau=h_{-1}=const-\frac{1}{8}\sinh^2{t}$.
This yields rather obvious

\medskip
{\bf Proposition 1.}
{\it Consider the round sphere
$(S^3,\bar{g}(\tau))$
satisfying the Ricci flow
$$
\frac{\partial}{\partial \tau}\bar{g}(\tau)=-2\bar{Ric}.
\eqno{(RF)}
$$
If the space coordinate~%
$t$
of the metric~(2)
is related to the flow time~%
$\tau$
by the function
$h_K$
then~
(2) is the metric of a~space of constant curvature~%
$K$ for $K\in\{-1,0,+1\}$.}

This claim is obvious.
Indeed,
consider the Ricci flow as a~procedure changing with time
the radius of the sphere
$S^3$.
Varying the rate of change of the radius,
we can ensure that
the space swept by the sphere
is of constant curvature.

Now we want to introduce some notations. Denote quotient $\frac{\partial_t\bar{g}}{\bar{g}}=2\frac{f'}{f}=2(\ln f)'$ as $\dot{f}$. And the quotient $\frac{\bar{Ric}}{\bar{g}}=\frac{4}{f^2}$ as $f_{Ric}$. Note that $(DF)$ will reduce to a single scalar equation because $g_{ij}=\delta_{ij}$ for orthonormal frame $\{\varepsilon^1,\varepsilon^2,\varepsilon^3\}$. Then the following theorem will be hold.

\medskip
{\bf Theorem 1.}
{\it If 
the round sphere
$(S^3,\bar{g}(t))$
changes under the action of the flow defined on the conformal factors
$$
\dot{f} =
\sqrt{f_{Ric}-4K},
$$
then the metric~%
$g$
is the metric of the space of constant curvature~%
$K$ for $K\in\{-1,0,+1\}$ if $f(0)=0$ and $f'(0)=1$.}

It is easy to see that system of equations from $(DF)$
is simple scalar equation
$$
2\frac{f'}{f}=\sqrt{\frac{4}{f^2}-4K}.
$$
It is clear that
for suitable initial data,
namely
$f(0)=0$
and
$f'(0)=1$,
this equation has the solutions
$f_{-1}$, $f_0$ and $f_{+1}$.
The associated metrics~%
$g$
are of the corresponding curvature~%
$K$.

One can check that
on the class of conformally round metrics
this flow and the standard Ricci flow
$(RF)$
with respect to time~%
$\tau$
related to the space coordinate~%
$t$
by
$h_K$
share the solutions describing spaces of constant curvature.

Remind that Dirac equation \cite{Fri}
$$
i\hbar \partial_t  = \sqrt{c^2\hbar^2\Delta + m^2c^4}
$$
is defined on the differential operators that act on some state functions $\psi$. It is more or less clear how to define the square root of the operator, and it is not clear how to define the square root of bilinear form $G$ such that $\sqrt{G}$ would be a bilinear form too. Also it is unclear what to do if the expression under the square root will be negative. 
Remind also that Ricci tensor can be regarded as Laplace operator of the metric in some very symmetric cases, and we think that $(DF)$ can be denoted as 'Dirac flow' although the resemblance is quite superficial.


\section[]{
A~metric with two functional parameters}

Consider an~embedding
$r_t:S^3\rightarrow \mathbb{R}^4$
inducing the metric
$$
g=dt^2+A_1^2(t)(e^1)^2+A_2^2(t)((e^2)^2+(e^3)^2)=dt^2+\bar{g}(t).\eqno{(4)}
$$
Now we define the orthonormal basis in slightly different way:
$$\varepsilon^0=dt,\quad \varepsilon^1=A_1(t)e^1,\quad \varepsilon^2=A_2(t)e^2,\quad \varepsilon^3=A_2(t)e^3.$$
Then the connection form is
$$
-\omega=\left(  \begin{array}{cccc}
                              0 & \frac{A_1'}{A_1}\varepsilon^1 & \frac{A_2'}{A_2}\varepsilon^2 & \frac{A_2'}{A_2}\varepsilon^3 \\
-\frac{A_1'}{A_1}\varepsilon^1 & 0 & \frac{-A_1}{A_2^2}\varepsilon^3 & \frac{A_1}{A_2^2}\varepsilon^2 \\
-\frac{A_2'}{A_2}\varepsilon^2 & \frac{A_1}{A_2^2}\varepsilon^3 & 0 & \frac{A_1^2-2A_2^2}{A_1A_2^2}\varepsilon^1 \\
-\frac{A_2'}{A_2}\varepsilon^3 & \frac{-A_1}{A_2^2}\varepsilon^2 & -\frac{A_1^2-2A_2^2}{A_1A_2^2}\varepsilon^1 & 0 \\
                            \end{array}
\right)
$$
and the curvature form is
$$
\Omega^0_1=\varepsilon^0\wedge\varepsilon^1[-\frac{A_1''}{A_1}]+ \varepsilon^2\wedge\varepsilon^3[-\frac{2A_1'}{A_2^2} +\frac{2A_1A_2'}{A_2^3}],
$$
$$
\Omega^0_2=\varepsilon^0\wedge\varepsilon^2[-\frac{A_2''}{A_2}]+ \varepsilon^3\wedge\varepsilon^1[\frac{A_1'}{A_2^2} -\frac{A_1A_2'}{A_2^3}],
$$
$$
\Omega^1_2=\varepsilon^0\wedge\varepsilon^3[\frac{A_1'}{A_2^2}-\frac{A_1A_2'}{A_2^3}]+ \varepsilon^1\wedge\varepsilon^2[\frac{A_1^2}{A_2^4}-\frac{A_1'A_2'}{A_1A_2}],
$$
$$
\Omega^2_3=\varepsilon^0\wedge\varepsilon^1[\frac{2A_1A_2'}{A_2^3}-\frac{2A_1'}{A_2^2}]+ \varepsilon^2\wedge\varepsilon^3[\frac{4}{A_2^2}-\frac{3A_1^2}{A_2^4}-\frac{A_2'^2}{A_2^2}],
$$
while the remaining components can be calculated similarly.
Observe that
for the metric~(4)
the curvature tensor is 
non-diagonal.
By a~diagonal tensor we consider a~curvature tensor
$R^i_{jkl}$
for which only the sectional curvatures
$R^i_{jij}$
can be nonzero.
Then the Ricci tensor of the metric~%
$\bar{g}$
has two nontrivial components:
$\bar{Ric}_{11}=\bar{R}^1_{212}+\bar{R}^1_{313}=4\frac{A_1^2}{A_2^4}$
and
$\bar{Ric}_{22}=\bar{Ric}_{33}=\bar{R}^1_{212}+\bar{R}^3_{232}=\frac{4}{A_2^2}(2-\frac{A_1^2}{A_2^2})$.

The Ricci flow reduces to the system of equations
$$
\left\{\begin{array}{ll}
\frac{A_1'}{A_1}\cdot\frac{dt}{d\tau}=-4\frac{A_1^2}{A_2^4}, \\
\frac{A_2'}{A_2}\cdot\frac{dt}{d\tau}=-\frac{4}{A_2^2}(2-\frac{A_1^2}{A_2^2}).
\end{array}\right.
\eqno{(5)}
$$
Assume firstly,
as in~(3),
that
$t\equiv \tau$.
Put
$\alpha=A_1^2$
and
$\beta=A_2^2$.
The system then becomes
$$
\left\{\begin{array}{ll}
\alpha'=-8\frac{\alpha^2}{\beta^2}, \\
\beta'+16=\frac{8\alpha}{\beta}.
\end{array}\right.
\eqno{(6)}
$$
This system has two obvious solutions.
The first corresponds to the case
$$
\alpha=\beta=8(t_0-t)
\eqno{(\mathrm{N})}
$$
already discussed.
The second solution
$$
\alpha=0,~\beta=16(t_0-t)
\eqno{(\mathrm{B})}
$$
corresponds to the collapse of the two-dimensional sphere
$S^2=S^3/S^1$
when the fiber
$S^1$
of the Hopf bundle is collapsed identically.
Essentially,
the second solution determines the metric~(4)
on a~three-dimensional space.

\medskip
{\bf Remark}.
We should elucidate our notation.
The solution~(N) corresponds to
the collapsing round three-dimensional sphere.
Gibbons and Hawking~\cite{GibHaw} call
solutions with this singularity ``nuts''
since in this case the whole three-dimensional sphere collapses.
A~different resolution of the cone singularity,
called a~``bolt'',
resembles the cylinder
$S^1\times\mathbb{R}$
with one coordinate lying on a~circle
$S^1$
of small radius
and the second coordinate is ''long'',
that is,
lying on the line~%
$\mathbb{R}$.
The second resolution of the singularity
amounts to a~circle
$S^1$
collapsing in the Hopf bundle with
the two-dimensional sphere of radius bounded below as a base.
In other words,
a~two-dimensional sphere plays the role of ''long'' coordinate,
while a~one-dimensional circle
$S^1$
(``bolt cross-section'')
collapses.
Our situation is more degenerate:
a~one-dimensional circle is collapsed identically,
$\alpha\equiv0$,
but in time the two-dimensional sphere collapses as well,
$\beta=16(t_0-t)$. Also we can note that in \cite{Baz} the similar situation appears although the space with singularity has dimension 8. The resolution of ''nut''-type is correspond to the space $\mathcal{M}_1$ and the ''bolt''-type singularity --- to the space $\mathcal{M}_2$.

\medskip
We can completely integrate~
(6).
From the second equation,
we find
$\alpha=\frac{1}{8}\beta(\beta'+16)$.
Inserting this into the first equation,
we obtain
$$
\beta\beta''+2\beta'^2+48\beta'+256=0.
$$
The general solution to this differential equation is 
$$
-\frac{1}{16}\beta-\frac{\sqrt{2}}{128}c_1\arctan(\frac{4\sqrt{2}\beta}{\sqrt{c_1^2-32\beta^2}})=t+c_2, \eqno{(7)}
$$
where
$c_1$
and
$c_2$
are constants of integration.
When the right-hand side of~(7) tends to zero from below,
the left-hand side tends to zero from above.
Furthermore,
$$
-\frac{c_1+\sqrt{c_1^2-32\beta^2}}{\sqrt{c_1^2-32\beta^2}}\cdot \frac{\beta'}{16}=1;
$$
hence,
for
$t\downarrow-c_2=t_0$
the derivatives satisfy
$\beta'\rightarrow -8$
and
$\alpha'=\frac{\beta'}{8}(\beta'+16)+\frac{\beta}{8}\beta''=\frac{\beta'}{8}(\beta'+16) +\frac{\beta}{8}(\frac{-1}{\beta})(256+48\beta'+2\beta'^2)\rightarrow -8$.
In other words,
we again recover a~well-known fact:
under the action of the Ricci flow
a~distorted three-dimensional sphere collapses to a~point
in finite time
(as a~round sphere of infinitesimal radius).
We can assert that
the solution~(B) is nonperturbative,
that is,
$$
\lim_{t\rightarrow t_0} \alpha'(t)=\lim_{t\rightarrow t_0} \beta'(t)=-8
$$
independently of how small
$\alpha(0)\neq 0$
is.
This agrees with general theory:
the Ricci flow makes curvature ``uniformly distributed''
over all points of the manifold and all tangent directions
provided that
the initial data is not too bad.
Therefore,
the sphere completely collapses
under the action of the Ricci flow;
the Ricci flow cannot decrease
just the one-dimensional fiber of the Hopf bundle,
keeping the radius of the two-dimensional base bounded below.
Thus,
we conclude that,
using the flows which differ ``insignificantly'' from the Ricci flow,
we have little chance to obtain a~description of metrics
with resolutions of conical singularity of ``bolt''-type.

Consider now the qualitative behavior of solutions.
If
$\beta(0)>\alpha(0)$
then at the initial moment of time
the right-hand sides of~
(5) are sufficiently small,
and~%
$\alpha$
behaves as a~constant,
while~%
$\beta$
decreases as
$-16t$.
When the radius
$A_2(t)$
of the sphere is sufficiently close to the radius
$A_1(t)$
of the circle,
they merge 
into one solution~%
$(N)$.
However,
if
$\alpha(0)>\beta(0)$
then,
since the right-hand side of the first equation
is a~negative number of large absolute value,
it follows that~%
$\alpha$
quite quickly becomes equal to
$\beta$,
and they merge again.

\medskip
\textbf{Remark}.
In computer simulations,
due to round-off errors
sometimes
(that is,
not in all simulations)
the solution to~
$(5)$
extends beyond the singularity time~%
$t_0$.
The function~%
$\alpha$
extends 
by zero,
while the function~%
$\beta$
becomes equal to
$16(t_0-t)$
and so goes into the negative region,
and the metric~%
$(4)$
ceases to be Riemannian;
moreover,
it ceases to be a~metric on a~four-dimensional manifold.
Nevertheless,
we can assert that
the solution~%
$(N)$
becomes the solution~%
$(B)$
by passing through the singularity.
It is not clear
whether we can find a~meaningful interpretation of this effect
or this is just an~artifact of simulations.

\medskip
This transformation of solution $(N)$ to $(B)$ has nothing to do with
the resolution of singularity
using the normalization of the Ricci flow.
Recall that
the normalized Ricci flow is 
$$
g'=-2Ric+\frac{2}{n}\mathfrak{R}g,
$$
where~%
$n$
is the dimension of the manifold
and~%
$\mathfrak{R}$
is the average scalar curvature.
In our case the normalized flow is
$$
\left\{\begin{array}{ll}
A_1'=-\frac{16}{3}\frac{A_1}{A_2^4}(A_1^2-A_2^2), \\
A_2'=\frac{8}{3}\cdot\frac{A_1^2-A_2^2}{A_2^3}
\end{array}\right.
$$
since the functions
$A_1$
and
$A_2$
depend only on time
and it is unnecessary to average the curvature over the sphere.
We can easily solve these equations;
more so when we remember that
the action of the normalized Ricci flow
preserves the volume
$\sqrt{det(\bar{g})}=A_1A_2^2$.
From the right-hand side of this system
it is not difficult to see also that,
as
$t\rightarrow \infty$,
both functions tend to the same constant.
Hence,
in infinite time we again obtain a~round sphere of constant radius.

We can define the decoupled 
Ricci flow~(5)
by replacing the left-hand side
while taking different derivatives with respect to time
in different directions:
$$
\left\{\begin{array}{ll}
\frac{A_1'}{A_1}\cdot r^{-1}=-4\frac{A_1^2}{A_2^4}, \\
\frac{A_2'}{A_2}\cdot s^{-1}=-\frac{4}{A_2^2}(2-\frac{A_1^2}{A_2^2}),
\end{array}\right.
\eqno{(5')}
$$
where
$r$~%
and~%
$s$
are two functions of~
$t$.
It obviously makes little sense to consider arbitrary
$r$~%
and~%
$s$;
therefore,
we try to find some geometrically meaningful particular cases.
Moreover, we can consider an~equation of even more general form:
$$
\mathfrak{D} g = R(g),
$$
where
$\mathfrak{D} g = (\frac{1}{r}\frac{\partial}{\partial t}g|_\mathcal{V}, \frac{1}{s}\frac{\partial}{\partial t}g|_\mathcal{H})$
is the derivative of the metric
determined by the natural Riemannian decomposition
into vertical and horizontal fibers of the space on which the considered metric is defined.
In many cases the horizontal subbundle
may also admit a~decomposition,
and the number of components of the derivative~%
$\mathfrak{D}$
may increase.

Recall that
a~metric is called anti-self-dual
whenever its connection form satisfies 
$$
\omega^i_j=-\frac{1}{2}\varepsilon_{ijkl}\omega^k_l,
$$
where~%
$\varepsilon$
is the Levi-Civita symbol.
The Ricci tensor of an~anti-self-dual metric,
as well as of a~self-dual one,
automatically vanishes.
In the 4-dimensional case
the anti-self-duality reduces to the pair of equations
$$
\omega^0_1=-\omega^2_3, \quad \omega^0_2=\omega^1_3
$$
or,
in our case,
$$
\frac{A_1'}{A_1}=-\frac{A_1^2-2A_2^2}{A_1A_2^2},\quad \frac{A_2'}{A_2}=\frac{A_1}{A_2^2}\eqno{(8)}.
$$

It is well-known that
these equations
are integrable;
this is how the classical Eguchi--Hanson metric \cite{EgHan}
$$
ds^2=[1-(a/r)^4]^{-1}dr^2 +r^2((e^2)^2+(e^3)^2)+r^2[1-(a/r)^4](e^1)^2
$$
was found.
This metric was the first metric with holonomy group
$SU(2)$
and an~explicit ansatz in elementary functions.
Recall that
as~%
$r$
tends to~%
$a$,
the Eguchi--Hanson metric has a~singularity of type ``bolt'',
while for sufficiently large~%
$R$
the set
$\{r=R\}$
is homeomorphic to
$\mathbb{R}P^3$. 

It is not difficult to observe that,
even though the right-hand sides of~
$(8)$
involve the components of the connection form,
they can be expressed in terms of the components of the Ricci tensor:
$$
-\frac{A_1^2-2A_2^2}{A_1A_2^2}=\frac{1}{2}\bar{Ric}_{22}(\bar{Ric}_{11})^{-1/2},
$$
$$
\frac{A_1}{A_2^2}=\frac{1}{2}(\bar{Ric}_{11})^{1/2}.
$$
Thus,
(8) is equivalent to the flow
$$
\frac{\partial}{\partial t}\bar{g}_{ij}=\frac{1}{2}{\sqrt{\mathrm{det}(\bar{Ric})}}{\bar{Ric}^{-1}_{ij}},
\eqno{(9)}
$$
where $\bar{Ric}^{-1}_{ij}$ are the components of the $3\times 3$-matrix of the Ricci tensor in orthonormal frame $\varepsilon^j$. We obtain a non-splitting
flow on the three-dimensional sphere,
that is,
a~flow which we can define on the total space of the Hopf bundle instead of defining it separately on the fiber
$S^1$
and base
$S^2$.
Unfortunately,
the resulting flow has an~extremely unpleasant right-hand side.
This confirms our suggestion that
a~metric with singularity of type ``bolt''
cannot be described by a~geometric flow with good right-hand side.
Observe also that
for
$A_1=A_2$
the function
$f(t)$
considered above
is a~solution,
while the corresponding metric~%
$g$
is the metric of a~flat space.
This follows immediately from the fact that
the metric of a~flat space is anti-self-dual.

\medskip
{\bf Theorem 2.}
{\it Consider the projective space
$\mathbb{R}P^3=S^3/{\mathbb Z}_2$
with the metric
$$
A_1^2(t)(e^1)^2+A_2^2(t)((e^2)^2+(e^3)^2)
$$
changing under the action of the flow~$(9)$.
If
$A_1(0)=0$,
$A_1'(0)=2$,
and
$A_2(0)=a$,
then the corresponding metric~%
$(4)$
is isometric to the Eguchi--Hanson metric.}

\medskip

As flow $(9)$ has such terrible right-hand side it would be natural to consider the evolution of some structures that respect the Hopf fibration. We remind that the 1-form $\psi$ on a $(2n+1)$-dimensional manifold called a contact form if $\psi \wedge (d\psi)^n \neq 0$. It is easy to check that
$$
\varepsilon^1\wedge \bar{d}\varepsilon^1=2\frac{A_1}{A_2^2}\varepsilon^1\wedge \varepsilon^2 \wedge \varepsilon^3,
$$
where $\bar{d}$ is a differential acting on $S^3$, that does not depends on $t$ and has no connection with Dolbeault theory. Also
$$
\varepsilon^2\wedge \bar{d}\varepsilon^2=2\frac{1}{A_1}\varepsilon^1\wedge \varepsilon^2 \wedge \varepsilon^3.
$$
If one defines the following flow on the contact structures
$$
(\ast \psi)'=\bar{d}\psi, \eqno{(10)}
$$
where $\ast$ is the Hodge operator (with respect to orthonormal frame $\{\varepsilon^1,\varepsilon^2,\varepsilon^3,\}$), then for $\psi=\varepsilon^1$ one will get exactly the second equation from$(8)$. And for $\psi=\varepsilon^2$ one will get
$$
(\frac{A_1'}{A_1}+\frac{A_2'}{A_2})\varepsilon^3\wedge \varepsilon^1 = \frac{2}{A_1}\varepsilon^3\wedge \varepsilon^1.
$$
That equation combined with one for $\psi\varepsilon^1$ will give the first equation from $(8)$. Obviously for $\psi=\varepsilon^3$ one will have exactly the same. So we have the following

\medskip
{\bf Theorem 3.}
{\it Consider the projective space
$\mathbb{R}P^3=S^3/{\mathbb Z}_2$
with the metric
$$
A_1^2(t)(e^1)^2+A_2^2(t)((e^2)^2+(e^3)^2)
$$
associated with 3 contact 1-forms $\varepsilon^1=A_1(t)e^1,\quad \varepsilon^2=A_2(t)e^2,\quad \varepsilon^3=A_2(t)e^3$, that are changed under the action of the flow $(10)$.
Then for appropriate initial data the metric
$$dt^2+A_1^2(t)(e^1)^2+A_2^2(t)((e^2)^2+(e^3)^2)$$
is isometric to the Eguchi--Hanson metric.}

\medskip
{\bf Remark} on connection with $G_2$.

We must mention that the flow $(10)$ was introduced by Hitchin (see for example \cite{Hitchin}) in order to obtain an evolution of the 3-form $\phi_0$ that defines the $G_2$-structure on 7-dimensional manifold. He proved that the solution $\phi_t$ of the system $(10)$ define a $G_2$-structure for sufficiently small times.

Also we observe that the equations $(8)$ with some trivial changes were written in \cite{BazMalk}. In that paper there was an attempt to construct a metric with $G_2$-holonomy on deformed cone over twistor space of 7-dimensional 3-Sasakian manifold.

\end{document}